\documentclass[11pt]{article} 
\usepackage{amsmath,amsthm} 
\usepackage{amssymb,latexsym}
\usepackage[T1]{fontenc}
\usepackage{authblk}
\usepackage{xcolor}
\usepackage{lipsum}
\usepackage{hyperref}
\usepackage{lipsum} 
\usepackage{lineno}
\modulolinenumbers[1]
\begin{document}
	\newtheorem{theorem}{Theorem}[section]
	\newtheorem{question}{Question}[section]
	\newtheorem{thm}[theorem]{Theorem}
	\newtheorem{lem}[theorem]{Lemma}
	\newtheorem{eg}[theorem]{Example}
	\newtheorem{prop}[theorem]{Proposition}
	\newtheorem{cor}[theorem]{Corollary}
	\newtheorem{rem}[theorem]{Remark}
	\newtheorem{deff}[theorem]{Definition}
	\numberwithin{equation}{section}
	\title{Structural Analysis of Commutative $S$-Reduced Rings} 
	
	\author[1]{Tushar Singh}
	\author[2]{ Shiv Datt Kumar }

	\affil[1, 2]{\small Department of Mathematics, Motilal Nehru National Institute of Technology Allahabad, Prayagraj 211004, India \vskip0.01in Emails: sjstusharsingh0019@gmail.com, tushar.2021rma11@mnnit.ac.in$^1$, sdt@mnnit.ac.in$^2$}
	
	\maketitle
	\hrule
	
	\begin{abstract}
		\noindent
	
	Let $R$ be a commutative ring with identity, $S \subseteq R$ be a multiplicative set. In this paper, we establish that the intersection of all $S$-prime ideals in an $S$-reduced ring is $S$-zero. Also, we show that an $S$-Artinian reduced ring is isomorphic to the finite direct product of fields. Furthermore, we provide an example of an $S$-reduced ring which is a uniformly-$S$-Armendariz ring (in short, $u$-$S$-Armendariz$)$ ring. Additionally, we prove that the class of uniformly-$S$-reduced rings (in short, $u$-$S$-reduced rings) belongs to the class of $u$-$S$-Armendariz rings. Among other results, we establish the relationship between $S$-reduced rings and $S$-strongly Hopfian rings. Finally, we prove the structure theorem for $S$-reduced rings.

	\end{abstract}
	\textbf{Keywords:} $S$-Reduced ring, $u$-$S$-Armendariz ring, $S$-integral domain, $S$-strongly Hopfian rings.\\
	\textbf{MSC(2020):} 12E20, 13G05, 16N40, 16U10.
	\hrule
	\section{Introduction}
	\noindent
	In this paper, $R$ denotes a commutative ring, and  $S\subseteq R$ is a multiplicative closed set. Commutative algebra and ring theory focus on the study of rings and their generalizations with the concept  of reduced rings being one of the most significant. Reduced rings generalize integral domains and serve as an intermediate framework making them instrumental in exploring properties between domains and general rings, preserved under generalization. In 2020, Pakin et al. \cite{ap20} introduced the notion of $S$-reduced rings as a generalization of reduced rings (see Definition \ref{newd}). The primary aim of this paper is to investigate various characterizations of $S$-reduced rings and to explore their relationships with certain special classes of rings, including   $u$-$S$-Armendariz rings, $S$-Hopfian rings, and $S$-PF rings. This study is motivated by results presented in \cite{ah20, ha21, nk00, ap20, ts23, ts24}. It is well known that, in a reduced ring, the intersection of all prime ideals is zero. However,  in an $S$-reduced ring, the intersection of all $S$-prime ideals  may not necessarily be zero (see Example \ref{eg1}). Recently, in 2024, Ersoy et al. \cite{ba24} introduced the concept of $S$-zero (see Definition \ref{zero}). Building on this idea, we establish that in an $S$-reduced ring, the intersection of all $S$-prime ideals  is precisely the $S$-zero (see Theorem \ref{szero}). Furthermore, in $2020$, Sevim et al. \cite{es19} extended the concept of Artinian rings to $S$-Artinian rings. A commutative ring $R$ is said to be \textit{$S$-Artinian} if for every descending chain of ideals:
	\begin{center}
		$I_{1}\supseteq I_{2}\supseteq\cdots\supseteq I_{n}\supseteq\cdots$
	\end{center}
	of $R$, there exist $s \in S$ and $k \in \mathbb{N}$ such that $s I_{k} \subseteq I_n$ for all $n \geq 1$. 
	
	In this section, we prove that an $S$-Artinian reduced ring is isomorphic to a finite direct product of fields (see Theorem \ref{art}). The concept of Armendariz rings was first introduced by Armendariz \cite{ea74} and later developed by Rege and Chhawchharia \cite{mb97}. These rings form an important class in algebra due to their ability to simplify and generalize polynomial ring theory, module theory, and other algebraic structures. Their importance lies in their distinctive properties and wide-ranging applications in both commutative and non-commutative ring theory.
	A ring \( R \) is called \textit{Armendariz} if whenever polynomials $f(x) = a_0 + a_1 x + \cdots + a_nx^n$, \quad $g(x) = b_0 + b_1 x + \cdots + b_m x^m\in R[x]$
	satisfy \( f(x)g(x) = 0 \), then \( a_i b_j = 0 \) for all \( 0 \leq i \leq n \) and \( 0 \leq j \leq m \). Due to its importance, the notion of a u-$S$-Armendariz ring, an excellent generalization of an Armendariz ring, was first presented by Kim et al. \cite{ki24} in $2024$ (see Definition \ref{armend}). In this work, we identify a specific class of \( 3 \times 3 \) upper triangular matrices that form a  $u$-\( S \)-Armendariz ring (see Example \ref{uarm}). Next, we establish that every $u$-$S$-reduced ring is a $u$-$S$-Armendariz ring (see Theorem~\ref{arm}).
	In \( 2021 \), Ahmed et al. \cite{ha21} introduced the concept of \( S \)-strongly Hopfian rings (see Definition \ref{strongly}). Subsequently, we demonstrate that every \( S \)-reduced ring is \( S \)-strongly Hopfian (see Proposition \ref{str}). Recently, in \( 2024 \), Ersoy et al. \cite{ba24} extended the notion of pure ideals to \( S \)-pure ideals (see Definition \ref{pure}). In this work, we show that the class of $S$-PF rings is a subclass of the class of $S$-reduced rings (see Example~\ref{pf}). Moreover, we conclude with a structure theorem for $S$-reduced rings (see Theorem~\ref{struct}).
	
	\section{$S$-Reduced Rings}
	We begin by introducing the concept of an \( S \)-reduced ring, with the goal of establishing  an \( S \)-analogue of well-known results associated with reduced rings.
	
	\begin{deff}$\cite{ap20}$\label{newd}
		A ring $R$ is said to be an $S$-\textit{reduced} if for any $r\in R$ such that $r^{n}=0$ for some $n \in\mathbb{N}$, then there exists $s\in S$ such that $sr= 0$.
	\end{deff}
	Clearly, every reduced ring is an $S$-reduced ring; however, the converse does not hold in general. To highlight this distinction, we now provide the following examples. Throughout the manuscript, we denote by $N(R)$ the set of all nilpotent elements of $R$
	
	\begin{eg}\label{r}
		Let $R=\mathbb{Z}_{24}$, and $S=\{\bar{1}, \bar{2}, \bar{4}, \bar{8}, \overline{16}\}$. We observe that $\bar{6}^{3}=\bar{0}$, but $\bar{6}\neq\bar{0}$. Therefore $R$ is not a reduced ring. Evidently, $N^*(R) = \{\bar{6}, \overline{12}, \overline{18} \}$ is the set of all non-zero nilpotent elements of $R$. Let $x\in N(R)$. Then we have $x^{3}=\bar{0}$. Take $s=\bar{4}\in S$. Then $sx=\bar{0}$, and hence $R$ is an $S$-reduced ring.
	\end{eg}
	\begin{eg}
		Let $R=\mathbb{Z}_{12}\times \mathbb{Z}_{12}\times \mathbb{Z}_{12}\times\ldots$ (countably infinite copies of $\mathbb{Z}_{12}$), and $S=\{\alpha^{n}\mid \hspace{0.1cm} where \hspace{0.1cm} \alpha=(\bar{2}, \bar{2}, \bar{2}, \ldots)\in R \hspace{0.1cm} and \hspace{0.1cm}n\in\mathbb{N}\cup\{0\}\}$. Let $\beta=(\bar{6}, \bar{6}, \bar{6}, \ldots)\in R$. We observe that $\beta^{2}=(\bar{0}, \bar{0}, \bar{0}, \ldots)$, but $\beta\neq(\bar{0}, \bar{0}, \bar{0}, \ldots)$. Thus $R$ is not a reduced ring. Clearly, the set $$ N(R)=\{(\bar{\alpha_{1}}, \bar{\alpha_{2}}, \bar{\alpha_{3}}, \ldots)|\  \alpha_{i}=0 \  \text{or} \  6 \ \text{for all}\ i \in \mathbb{N}\}$$ is the set of all nilpotent elements of $R$. Let $\gamma\in N(R)$. Then $\gamma^{n}=(\bar{0}, \bar{0}, \bar{0}, \ldots)$ for some $n\in \mathbf{N}$. Take $\alpha=(\bar{2}, \bar{2}, \bar{2}, \ldots)\in S$. Consequently, $s\gamma=(\bar{0}, \bar{0}, \bar{0}, \ldots)$ and hence $R$ is an $S$-reduced ring. 
	\end{eg}
	\begin{rem}
		If $S^{-1}R$ is an integral domain, then $R$ is an $S$-reduced ring. For this, let $r\in R$ be such that $r^{n}=0$ for some $n\in\mathbb{N}$. Then we have $\dfrac{r^{n}}{s^{n}} =\dfrac{0}{s^{n}}$ for some $s\in S$, $\dfrac{r}{s} =\dfrac{0}{s}$ since $S^{-1}R$ is an integral domain. This implies that there exists $u\in S$ such that $usr=0$, as desired.
	\end{rem}
	
	\begin{deff}\cite[Definition 2.1]{ut25}
		Let $S\subseteq R$ be a multiplicative set and $I$ be an ideal of $R$. Then the $S$-radical of $I$ is defined by
		$$\sqrt[S]{I} = \{a \in R \mid sa^n \in I \text{ for some } s \in S \text{ and } n \in \mathbb{N}\}. 
		$$ Further,  $I$ is said to be an $S$-radical ideal if $\sqrt[S]{I} =I$.
	\end{deff}
\noindent
Now $Z(R)$ denotes the set of zero divisors of a ring $R$.
	\begin{prop}\label{sr}
		Let $I$ be an ideal of $R$ with $Z(R/I)\cap \bar{S}=\emptyset$, where $\bar{S}=\{s+I\mid s\in S\}$ be a multiplicative closed set of $R/I$. Then  $I$ is an $S$-radical ideal if and only if $R/I$ is an $\bar{S}$-reduced ring.
	\end{prop}
	\begin{proof}
		If $I$ is $S$-radical, then it is radical and so $R/I$ is reduced, which means it is $\bar{S}$-reduced ring. 
		Conversely, let $R/I$ be an $\bar{S}$-reduced ring, and $y\in \sqrt[S]{I}$. Then there exists $s\in S$ such that $sy^{n}\in I$, $(sy)^{n}\in I$. Consequently, $(sy+I)^{n}=I$; there exists $t+I\in \bar{S}$ such that $(t+I)(sy+I)=I$ since $R/I$ is an $\bar{S}$-reduced ring. It follows that $(st+I)(y+I)=I$, then $y\in I$ since $Z(R/I)\cap \bar{S}=\emptyset$, as desired.
	\end{proof}

	
	\begin{question}
		Is the assumption \textquotedblleft $S\cap Z(R/I)=\emptyset$\textquotedblright \hspace{0.1cm} necessary for
		the Proposition \ref{sr}?
	\end{question}
	Yes, this is a necessary condition. We provide an example in which $S\cap Z(R/I)\neq \emptyset$, and Proposition \ref{sr} is not valid.
	\begin{eg}
		Let $R=\mathbb{Z}$, $S=\{2^{n}\mid n\in\mathbb{N}\cup \{0\}\}$, and $I=12\mathbb{Z}$. Then $\overline{R}=R/I=\mathbb{Z}_{12}$, $\bar{S}=\{s+12\mathbb{Z}\mid s\in S\}=\{\bar{1}, \bar{2}, \bar{4}, \bar{8}\}$. Evidently, $Z(\overline{R})\cap \bar{S}\neq\emptyset$. Since $\bar{6}$ is the only non-zero nilpotent element of $\overline{R}$. Take $s=\bar{2}\in \overline{S}$. Then $s\bar{6}=\bar{0}$, and hence $\overline{R}$ is an $S$-reduced. Since $3\notin I$ but $\bar{4}\cdot 3\in I$, it follows that $I$ is not an $S$-radical ideal of $R$.
	\end{eg}
	
	\begin{deff} $\cite{ba24}$\label{zero}
		Let $y\in R$. Then we have the following:
		\begin{enumerate}
			\item  If $sy^{n}=0$ for some $s\in S$ and $n\in \mathbb{N}$, then $y$ is called $S$-nilpotent.
			\item  If $sy = 0$ for some $s\in S$, then $y$ is called $S$-zero.  
		\end{enumerate}
	\end{deff}
\begin{deff}
An ideal \( I \) of a ring \( R \) is called \emph{\( S \)-zero} if every element of \( I \) is \( S \)-zero; that is, for each \( a \in I \), there exists \( s \in S \) such that \( sa = 0 \).

\end{deff}
	
	\begin{prop}
		Let $R$ be an $S$-reduced ring and  $I$, $J$ be ideals of $R$. Then $I\cap J$ is an $S$-zero ideal if and only if $IJ$ is an $S$-zero ideal.
	\end{prop}
	\begin{proof}
		Assume that $I\cap J$ is an $S$-zero ideal of $R$. Let $a\in IJ\subseteq I\cap J$. This implies that there exists $s\in S$ such that $sa=0$. Thus $IJ$ is an $S$-zero ideal of $R$. Conversely, let $a\in I\cap J$. Then $a\in I$ and $a\in J$, $a^{2}\in IJ$. Consequently, there exists $s'\in S$ such that $s'a^{2}=0$ since $IJ$ is an $S$-zero ideal of $R$. Since $R$ is an $S$-reduced ring, there exists $s''\in S$ such that $s's''a=0$. Therefore $I\cap J$ is an $S$-zero ideal of $R$.
	\end{proof}
	
	\begin{deff}$\cite{ed20}$\label{sintegral} 
		A ring $R$ is called an $S$-\textit{integral domain} if there exists $s \in S$ such that for any $a, b \in R$, if $ab = 0$, then  either $sa = 0$ or $sb = 0$.
	\end{deff}
	
	Note that every integral domain is $S$-integral domain; however, the converse does not hold in general (see Example \ref{s})).
	\begin{eg}\label{s}
		Consider $R=\mathbb{Z}\times \mathbb{Q}$, ~$S=\mathbb{Z}\times \mathbb{Q}^{*}$, where $\mathbb{Q}^{*}$ is a non-zero rational number. Evidently, $R$ is not an integral domain. Let $a, ~b\in R$ such that $a=(x, u)$ and $b=(y, v)$, where $x, y\in \mathbb{Z}$ and $u, v\in\mathbb{Q}$. Now, if $a\cdot b=(xy, uv)=(0, 0)$, then it follows that $xy=0$ and $uv=0$. This implies that either $x=0$ or $y=0$ and either $u=0$ or $v=0$ since $\mathbb{Z}$ and $\mathbb{Q}$ are integral domains. Therefore, there are two possibilities for $a$ and $b$, if $a=(x, 0)$, then $b=(0, v)$. Further, if $a=(0, u)$, then $b=(y, 0)$. Take $s=(0,1)\in S$. Then either $sa = 0$ or $sb = 0$. Thus $R$ is an $S$-integral domain.
	\end{eg}
	\noindent
	\textbf{Note:} It is a trivial fact that all $S$-integral domains are $S$-reduced rings.
	
	\noindent
	\begin{deff}\cite{fm69} Let $I$ and $J$ be two ideals of $R$. 
		\begin{enumerate}
				\item  Then the \textit{colon ideal} of $I$ by $J$ is defined as $(I : J) = \{ r \in R \mid rJ \subseteq I \}$.

			\item In particular, if $x \in R$, then the colon ideal of $I$ by $x$ is defined as $(I : x) = \{ r \in R \mid rx \in I \}$.
		
		\end{enumerate}
	\end{deff}
	
	\noindent
	Following from \cite{ah20}, an ideal $P$ disjoint from  $S$ of a ring $R$ is said to be an $S$-\textit{prime }if there exists an $s\in S$ such that for all $a,b\in R$ with $ab\in P$, we have $sa\in P$ or $sb\in P$ or equivalently, $P$ is $S$-\textit{prime} if and only if there exists $s \in S$, such that for all $I, J$ two ideals of $R$, if $IJ \subseteq P$, then $sI \subseteq P$ or $sJ \subseteq P$.
	Clearly, every prime ideal is an $S$-prime, but the converse is not true in general (see \cite[Example 1(3)]{ ah20}). We know that the intersection of all prime ideals in a reduced ring is zero. Then we have the following question: 
	\begin{question}
		Is the intersection of all $S$-prime ideals of an $S$-reduced ring equal to zero ?
	\end{question}
The answer to this question is negative. We present an example showing that the intersection of all $S$-prime ideals of an $S$-reduced ring need not be the zero ideal. Before giving the example, we include the following remark to provide some necessary context:
\begin{rem}
Recall from \cite{em60} that if \( D \) is an integral domain with quotient field \( \mathcal{F} \), and \( M \) is a \( D \)-module. Then \( M \) is called a \textit{divisible} \( D \)-module, if \( rM = M \) for every \( 0 \neq r \in D \). Furthermore, following \cite[Proposition~1]{ah20}, an ideal \( P\) of a ring \( R \) that is disjoint from \( S \) is an \( S \)-\textit{prime ideal} of \( R \) if and only if there exists \( s \in S \) such that \( (P : s) \) is a prime ideal of \( R \).
\end{rem}
 Following from \cite{ad09},  let $M$ be an $R$-module. The \emph{idealization} of the $R$-module $M$, denoted by $R(+)M$, is defined as $R(+)M = \{(r, m) \mid r \in R,\ m \in M\}$, which forms a commutative ring with componentwise addition and multiplication given by $(\alpha_1, m_1)(\alpha_2, m_2) = (\alpha_1 \alpha_2,\, \alpha_1 m_2 + \alpha_2 m_1)$, for all $\alpha_1, \alpha_2 \in R$ and $m_1, m_2 \in M$. It is straightforward to verify that $S(+)N = \{(s, n) \mid s \in S,\ n \in N\}$, where $N$ is a submodule of $M$, forms a multiplicative set in $R(+)M$.

\begin{eg}\label{eg1}
	Let $D$ be an integral domain, which is not a field and $\mathcal{F}$ denotes its quotient field. Define $R=D(+)(\mathcal{F}/D)$ and $S$ the set of all non-nilpotent elements of R. Then $S=\{(x, (a/b)+D)\mid x\not=0,a,b\in D, b\not=0\}$ is a multiplicative subset of $R$. Let $r\in R$ with $r^n=0$. Then $r=(0,a/b+D)$ for some $a,b\in D, b\not=0$. Put $s=(b,0+D)\in S$. Then $(b,0+D)r=(b,0+D)(0,a/b+D)=(0, D)=0_{(+)}$, where $0_{(+)}$ denotes the zero elements of $R$, and so $R$ is $S$-reduced. Now, we prove that $R$ is not an $S$-integral domain. On the contrary, assume that $R$ is an $S$-integral domain. Then $(0_{(+)})$ is an $S$-prime ideal of $R$, and so $\left((0_{(+)}):(x,a/b+D)\right)$ is a prime ideal of $R$, by \cite[Proposition 1]{ah20}. This implies that $\left((0_{(+)}):(x,a/b+D)\right)=p(+)(\mathcal{F}/D)$ for some prime ideal $p$ of $D$. Since $x\not=0$, we have $p=(0)$. Hence $\left((0_{(+)}):(x,a/b+D)\right)=(0)(+)(\mathcal{F}/D)$. Note that $(x, a/b+D)\left((0)(+)(\mathcal{F}/D)\right)=(0)(+)(\mathcal{F}/D)\not=(0_{(+)})$ as $\mathcal{F}/D$ is divisible $(x\mathcal{F}/D=\mathcal{F}/D)$, which is a contradiction. Hence $R$ is not an $S$-integral domain. Next, we show that the intersection of all $S$-prime ideals of $R$ is non-zero. Let $P$ be an $S$-prime ideal of $R$. Since $P\cap S=\emptyset$, we have $P= (0)(+)N$ for some $D$-submodule $N$ of $\mathcal{F}/D$. Then $\left((0)(+)N:(x,a/b+D)\right)=p(+)(\mathcal{F}/D)$ for some prime ideal $p$ of $D$. As above, we have $p=(0)$. So $\left( (0)(+)N:(x,a/b+D)\right)=(0)(+)(\mathcal{F}/D)$. Since $\mathcal{F}/D$ is divisible, $(x, a/b+D)\left((0)(+)\mathcal{F}/D\right)=(0)(+)(\mathcal{F}/D)$. Hence $N=\mathcal{F}/D$. Consequently, the only $S$-prime ideal of $R$ is $(0)(+)(\mathcal{F}/D)$. Hence, the intersection of all $S$-prime ideals of $R$ is $(0)(+)(\mathcal{F}/D)$, which is non-zero.
\end{eg}

We are now prepared to extend this property to the $S$-version. Here, Spec$_S(R)$ denotes the collection of all $S$-prime ideals of $R$.
	\begin{theorem}\label{szero}
		The intersection of all $S$-prime ideals of an $S$-reduced ring is $S$-zero provided $0\notin S$.
	\end{theorem}
	\begin{proof}
		On the contrary suppose $\bigcap_{P\in Spec_{S}(R)} P$ is not an $S$-zero.
		Then there exists $x\in\bigcap_{P\in Spec_{S}(R)} P$ such that $sx\neq 0$ for all $s\in S$. This implies that $x^{n}\neq 0$ for all $n\in\mathbb{N}$ since $R$ is an $S$-reduced ring. Now, define a set $$E=\{sx^{n}\mid s\in S, n\in\mathbb{N}\cup\{0\}\}.$$ Again, consider a collection of ideals defined with the help of set $E$ as follows: $$E'=\{I\subseteq R\mid I\cap E=\emptyset\}.$$ Clearly, $E'\neq \emptyset$, as $(0)\in E'$. By Zorn's lemma, $E'$ has a maximal element, say $P'$. To prove that $P'$ is an $S$-prime ideal of $R$, first we prove that $P'\cap S=\emptyset$. For this, we choose $n=0$, and abotain $sx^n=s$, so all elements of $S$ are included in $E$. Thus $P'$ is disjoint from $S$ because it belongs to $E'$. Now, we claim $(P':t)\cap E=\emptyset$ for all $t\in S$. Contrary, suppose that there exists $s'\in S$ such that $(P':s')\cap E\neq\emptyset$. Consequently, there exists $v\in E$ such that $v\in (P':t)$. It follows that $stx^{l}\in P'$ for some $l\in\mathbb{N}$, a contradiction, as $P'\cap E=\emptyset$. Thus $(P':t)\cap E=\emptyset$ for all $t\in S$. Next, we claim that $P'$ is $S$-prime ideal with respect to $u\in S$. Then, by \cite[Proposition 1]{ah20}, it is sufficient to show that $(P':u)$ is a prime ideal of $R$. Suppose $(P':u)$ is not a prime ideal. Then there exist $a, b\in R$ such that $ab\in(P':u)$ but neither $a\in (P':u)$ nor $b\in(P':u)$. This implies that $ua\notin P'$ and $ub\notin P'$, and so $P'\subseteq P'+Rua$ and $P'\subseteq P'+Rub$. Then $(P'+Rua)\cap E\neq \emptyset$ and $(P'+Rub)\cap E\neq \emptyset$ since $P'$ is a maximal element of $E'$. Consequently, there exist $v, w\in E$ and $l, m\in \mathbb{N}$ such that $v=sx^{l}=r_1ua+c$ and $w=r_{2}ub+d$ for some $r_{1}, r_{2}\in R$ and $c, d\in P'$. Now, we have $s^{2}x^{l+m}=r_{1}r_{2}u^{2}ab+u(r_1ad+r_2ubc)+cd\in (P':u)$, a contradiction, as $(P':t)\cap E=\emptyset$ for all $t\in S$. Thus $(P':u)$ is a prime ideal of $R$, and hence $P'$ is $S$-prime with respect to $u\in S$, by \cite[Proposition 1]{ah20}. Thus $x\in P'$; $sx^{n}\in P'$, a contradiction, as $P'\cap E=\emptyset$. Hence $\bigcap_{P\in Spec_{S}(R)} P$ is an $S$-zero ideal.
	\end{proof}
	\begin{deff} \label{t} Let $I$ be an ideal of a ring $R$. Then we have the following: 
		\begin{enumerate}
			\item  If for any
			$a\in I$, there exists $s\in S$ such that $sa^n = 0$  for some $n\in\mathbb{N}$,  then $I$ is called an $S$-\textit{nil ideal}.
			\item The set $Nil_{S}(R)=\{a\in R\mid sa^n = 0 ~for ~some ~n\in \mathbb{N} ~and ~s\in S\}$.	Obviously,  $Nil_{S}(R)$ is an ideal of $R$, which is called the $S$-Nil radical of $R$.
		\end{enumerate}
	\end{deff}

	\begin{rem}\label{pre}
		Let $P\in Spec_{S}(R)$. Then, by \cite[Proposition 2]{ed20}, there exists $s\in S$ such that $(P: s')\subseteq (P: s)$ for all $s'\in S$ and $(P: s)$ is a prime ideal. From now on, we denote this $s\in S$ for $P\in Spec_{S}(R)$ by $s_{P}$.
	\end{rem}
	\begin{cor}
		Let $R$ be a ring, $S\subseteq R$ be a multiplicative set. Then we have the following:
		\begin{enumerate}
			\item If $P$ is an $S$-prime ideal of $R$, then $Nil_{S}(R)\subseteq (P:s_{P})$.
			\item If $R$ is an $S$-reduced ring, then $Nil_{S}(R)$ is $S$-zero.
		\end{enumerate}
	\end{cor}
	\begin{proof}
		\leavevmode
		\begin{enumerate}
			\item Let $a\in Nil_{S}(R)$. Then there exists $s'\in S$ such that $s'a^{n}=0$, $(s'a)^{n}\in P$. Since $P$ is an $S$-prime ideal of $R$, there exists $s\in S$ such that $ss'a\in P$, and hence $Nil_{S}(R)\subseteq (P:s_{p})$, by Remark \ref{pre}. 
			\item Immediate from the Definition \ref{t}(2).
		\end{enumerate}
	\end{proof}

	\begin{rem}\label{red}
		If $R$ is an $S$-reduced ring, then $S^{-1}R$ is a reduced ring.
	\end{rem}
	\begin{proof}
		Let $0\neq x=\frac{a}{s}\in S^{-1}R$, where $a\in R$ and $s\in S$ be such that $x^{n}=0$ for some $n\in \mathbb{N}$. Then $\frac{a^{n}}{s^{n}}=\frac{0}{1}$. Consequently, there exists $u\in S$ such that $ua^{n}=0$, $(ua)^{n}=0$. Since $R$ is $S$-reduced, there exists $u'\in S$ such that $uu'a=0$. Then $x=\frac{uu'a}{uu's}=\frac{0}{1}$. Therefore $S^{-1}R$ is a reduced ring.
	\end{proof}
The converse of Remark~\ref{red} need not hold. An illustrative example is given below.

	\begin{eg}\label{nsint}
		\noindent
	Let $\mathbb{E}=\bigoplus_{p\in\mathcal{P}}\mathbb{Z}/p\mathbb{Z}$, where $\mathcal{P}$ is the set of all prime numbers in $\mathbb{Z}$. Define $R=\mathbb{Z}(+)\mathbb{E}$ and $S=\mathbb{Z}\setminus\{0\}(+)(0)$. According to \cite[Example 3.12]{ki24}, the localization 
	$S^{-1}R=\mathbb{Q}$, which is a field and hence a reduced ring. Let $\alpha=(0, 1)\in R$. Then $\alpha^{2}=0$ but $s\alpha\neq 0$ for all $s\in S$, for this, if there exists $s\in S$ such that $s\alpha=0$. Evidently, $s=(n, 0)\in S$, where $n\in\mathbb{Z}\setminus\{0\}$. It follows that $s\alpha=(0, n)=0$, and so $n=0$, which is not possible since $n\in\mathbb{Z}\setminus\{0\}$. Therefore $R$ is not an $S$-reduced ring.
	\end{eg}

Following  from \cite{ad02}, an ideal $I$ of a ring $R$ is called $S$-\textit{finite}, where $S\subseteq R$ be a multiplicative set, if $sI \subseteq J \subseteq I$ for some finitely generated ideal $J$ of $R$ and some $s \in S$. If each ideal of $R$ is $S$-finite, then  $R$ is called $S$-\textit{Noetherian ring}.

	\begin{theorem}\label{art}
		Let $R$ be an $S$-Noetherian $S$-reduced ring with the property that every nonzero element of $S^{-1}R$ is a zero divisor or invertible. Then $S^{-1}R$ is an Artinian ring.
	\end{theorem}
	\begin{proof}
		Note that $S^{-1}R$ is a Noetherian ring, by \cite[Proposition 2(f)]{ad02}. Then $S^{-1}R$ has a finite number of minimal prime ideals, say $\{S^{-1}P_{1},\ldots, S^{-1}P_{n}\}$, where $P_{i}$ $(1\leq i\leq n)$ is an $S$-prime ideal of $R$, by \cite[Corollary 4]{ah20}. By Remark \ref{red}, $S^{-1}R$ is a reduced ring, $\{0\}=\bigcap\limits_{i=1}^{n}S^{-1}P_{i}$ since in reduced ring intersection of all minimal prime ideals is zero. Suppose $S^{-1}P$ is a prime ideal of $S^{-1}R$. Let $0\neq x\in S^{-1}P$. Then there exists $0\neq y\in S^{-1}R$ such that $xy=0$. This implies that $xy\in S^{-1}P_{i}$ for all $i=1,\ldots, n$. If $x\notin S^{-1}P_{i}$ for all $i=1,\ldots, n$, then $y\in\bigcap\limits_{i=1}^{n}S^{-1}P_{i}$, a contradiction as $y\neq 0$. Thus $x\in S^{-1}P_{i}$ for some $i=1,\ldots, n$. It follows that $S^{-1}P\subseteq \bigcup\limits_{i=1}^{n} S^{-1}P_{i}$. By \cite[Proposition 1.11]{fm69}, we get $S^{-1}P\subseteq S^{-1}P_{i}$ for some $i$, $S^{-1}P=S^{-1}P_{i}$. Hence $S^{-1}R$ has only finitely many prime ideals of the form $\{S^{-1}P_{i}\}$ and  $dim(S^{-1}R)=0$. Thus $S^{-1}R$ is an Artinian ring.
	\end{proof}
	Recall from \cite{es19} that a ring $R$ is said to be \textit{$S$-Artinian} if for every descending chain of ideals $I_{1}\supseteq I_{2}\supseteq\cdots\supseteq I_{n}\supseteq\cdots$ of $R$, there exist $s \in S$ and $k \in \mathbb{N}$ such that $s I_{k} \subseteq I_{n}$ for all $n \geq k$.
	We denote by $Nil(R)$ the intersection of all prime ideals of $R$, called the nilradical. It is well-known that in an Artinian ring, the nilradical is nilpotent (see \cite[Proposition~8.4]{fm69}). Ansari and Sharma \cite{ab} extended this result to the setting of $S$-Artinian rings as follows:
%

	\begin{lem}$\cite{ab}$\label{nil}
		Let $R$ be a ring, $S\subseteq R$ be a multiplicative closed set which contains no zero-divisors of $R$. Then we have the following:
		\begin{enumerate}
			\item $Nil(R)$= Intersection of all prime ideals of $R$ which are disjoint from $S$.
	
			\item If $R$ is an $S$-Artinian ring, then $Nil(R)$ is nilpotent.
		\end{enumerate}
	\end{lem}
	
	\begin{theorem}
		An $S$-Artinian reduced ring is isomorphic to the finite direct product of fields provided $Z(R)\cap S=\emptyset$.
	\end{theorem}
	\begin{proof}
		Assume that $R$ is an $S$-Artinian and reduced ring. According to \cite[Theorem 2.4]{es20}, in the $S$-Artinian ring there are only finitely many prime ideals disjoint from $S$, say, $P_{1}, P_{2}, \ldots, P_{n}$, then we have $Nil(R)=\bigcap\limits_{i=1}^{n}P_{i}$, by Lemma \ref{nil}(1). Since $Nil(R)$ is nilpotent in an $S$-Artinian ring, there exists $k\in\mathbb{N}$ such that $(\bigcap\limits_{i=1}^{n}P_{i})^{k}=0$, by Lemma \ref{nil}(2). Since $\prod_{i=1}^{n}P_{i}\subseteq \bigcap\limits_{i=1}^{n}P_{i}$, then $(\prod_{i=1}^{n}P_{i})^{k}=0$. This implies that $(\prod_{i=1}^{n}P_{i})=0$ since $R$ is a reduced ring. Also, by $\cite[Proposition ~2.5]{es20}$, each $P_{i}$ $(1\leq i\leq n)$ is minimal and also maximal ideal of $R$, by \cite[Proposition 8.1]{fm69}. Then by Chinese Remainder Theorem, we have:
		$$R=R/(\prod_{i=1}^{n}P_{i})\cong R_1\times R_{2}\times \cdots \times R_{n}, $$ where $R_{i}=R/P_{i}$ is a field for all $i=1, 2, \ldots, n$. 
	\end{proof}
\noindent
	\textbf{Note:} It is well-known that a ring $R$ is reduced if and only if the polynomial ring $R[X]$ is reduced. We now extend this result to the notion of $S$-reduced rings:
	\begin{prop}
		A ring $R$ is $S$-reduced ring if and only if  $R[X]$ is an $S$-reduced ring.
	\end{prop}
	\begin{proof}
		Suppose that $R$ is an $S$-reduced ring. Let $f(X)=\sum_{i=1}^{n}a_{i}X^{i}\in R[X]$, where $a_{i}\in R$ for $(1\leq i\leq n)$ such that $(f(X))^{k}=0$ for some $k\in\mathbb{N}$. Then there exist $s_{i}\in S$ for each $a_{i}$ such that $s_{i}a_{i}=0$. Take $s=s_{1}s_2\ldots s_{n}\in S$. Then $sa_{i}=0$. Thus $R[X]$ is $S$-reduced. Conversely, suppose that $R[X]$ is an $S$-reduced ring. Let $a\in R$ such that $a^n = 0$ for some $n \in \mathbb{N}$. Consider the polynomial $f(X) = a \in R[X]$. Since $f(X)^n = a^n = 0$ in $R[X]$ and $R[X]$ is $S$-reduced, there exists $s\in S$ such that $sa= 0$. Hence, $R$ is an $S$-reduced ring.
	\end{proof}
	\section{Uniformly-$S$-Armendariz Rings}
In 1974, E. Armendariz introduced the class of rings now known as \emph{Armendariz rings}, and showed that every reduced ring is Armendariz (see \cite[Lemma~1]{ea74}). More recently, in 2024, Kim et al.\ extended this notion to the framework of multiplicative sets by introducing the concept of \emph{$u$-$S$-Armendariz rings}. Motivated by this, we begin this section with an example of an $S$-reduced ring that is also a $u$-$S$-Armendariz ring.
%
	
	\begin{deff}$\cite{ki24}$\label{armend}
		A ring $R$ is called an $u$-$S$-\textit{Armendariz} if for any polynomials $f =\sum_{i=1}^{m} a_iX^i$ and $g =\sum_{j=1}^{n} b_jX^j$ in $R[X]$ satisfying $fg = 0$, then there exists an element $s\in S$ such that $sa_ib_j = 0$ for all indices $i$ and $j$.
	\end{deff}
	\begin{eg}\label{uarm}
		Let $R$ be an $S$-reduced ring. Then 
		\[E= \left\{\begin{bmatrix}
		a & b & c \\
		0 & a & d\\
		0 & 0 & a
		\end{bmatrix}\mid a, b,c, d\in R \right\}	
		\] is an u-$S'$-Armendariz ring, where \[S'= \left\{\begin{bmatrix}
		s& s& s \\
		0 & s & s \\
		0 & 0 & s
		\end{bmatrix}\mid s\in S \right\}	
		\] is a multiplicative closed subset of $E$.
	\end{eg}
	\begin{proof}
		\noindent
		We employ the method in the proof of \cite[Proposition 2]{nk00}. First, we notice that for 
		\[
		\begin{bmatrix}
		a_{1} & b_{1} & c_{1}  \\
		0 & a_{1} & d_{1} \\
		0 & 0 & a_{1} 
		\end{bmatrix}, \quad
		\begin{bmatrix}
		a_{2} & b_{2} & c_{2}  \\
		0 & a_{2} & d_{2} \\
		0 & 0 & a_{2} 
		\end{bmatrix} \in E,
		\]
		we can denote their addition and multiplication as: $(a_1, b_1, c_1, d_1)+(a_2, b_2, c_2, d_2)= (a_1 + a_2, ~b_1 + b_2, ~c_1 + c_2, ~d_1 + d_2)$ and $(a_1, b_1, c_1, d_1)(a_2, b_2, c_2, d_2)= (a_1 a_2, ~a_1 b_2 + b_1 a_2, ~a_1 c_2 + b_1 d_2 + c_1 a_2, ~a_1 d_2 + d_1 a_2)$ respectively. So every polynomial in $E[x]$ can be expressed in the form  $\big(p_0(x), ~p_1(x), ~p_2(x), ~p_3(x)\big)$ for some $p_i(x)\in R[x]$.
		
		Let $f(x) = \big(f_0(x), ~f_1(x), ~f_2(x), ~f_3(x)\big)$ and $g(x) = \big(g_0(x), ~g_1(x), ~g_2(x), ~g_3(x)\big)$ be elements of $E[x]$ such that $f(x)g(x)= 0$. It follows that $f(x)g(x) = \big(f_0(x)g_0(x), ~f_0(x)g_1(x)+f_1(x)g_0(x), ~f_0(x)g_2(x)+ f_1(x)g_1(x)+f_2(x)g_0(x), ~f_0(x)g_3(x)+f_1(x)g_2(x)\big)= 0$.
		Thus we have the system of equations:
		\[
		\begin{aligned}
		\quad f_0(x) g_0(x) = 0, ~~~~~~~~~~~(A)\\
		\quad f_1(x) g_0(x) + f_0(x) g_1(x) = 0,~~~~~~~~~~~(B) \\
		\quad f_2(x) g_0(x) + f_1(x) g_3(x) + f_0(x) g_2(x) = 0, ~~~~~~~~~~~(C)\\
		\quad f_3(x) g_0(x) + f_0(x) g_3(x) = 0. ~~~~~~~~~~~(D)
		\end{aligned}
		\] 
		Now, if we multiply  equation \((B)\) by \(f_0(x)\), then $f_1(x)f_0(x)g_0(x) + f^{2}_0(x) g_1(x)= 0$, $(f_0(x) g_1(x))^{2}= 0$. Using the fact that \( R[x]\) is an $S$-reduced, there exists $s_{1}\in S$ such that $s_{1}f_0(x) g_1(x)= 0$, and hence $s_{1}f_1(x) g_0(x)= 0$. Similarly, if we multiply  equation \((D)\) by \(f_0(x)\), then $f_3(x)f_0(x)g_0(x) + f^{2}_0(x)g_3(x)= 0$, $(f_0(x)g_3(x))^{2}= 0$. Then there exists $s_{2}\in S$ such that $s_{2}(f_0(x)g_3(x))=0$, and hence $s_{2}(f_3(x)g_0(x))=0$. Now, if we multiply  equation \((C)\) by $s_{2}$\(f_0(x)\), then $s_2f_2(x)f_0(x)g_0(x) + f_1(x)s_2f_0(x)g_3(x) + s_2f^{2}_0(x)g_2(x) = 0$, $(s_2f_0(x)g_2(x))^{2}= 0$. Then there exists $s_3\in S$ such that $s_2s_{3}(f_0(x)g_2(x))= 0$, and hence  equation \((C)\) becomes  $$s_2s_{3}(f_2(x) g_0(x)+f_1(x) g_3(x))= 0. ~~~~~~~~~~~(D')$$ Again multiply equation \((D')\) by $s_{1}$\( f_1(x) \), then $s_2s_{3}(f_2(x) s_{1}f_{1}(x)g_0(x)+s_1f^{2}_1(x) g_3(x))= 0$, $(s_1s_2s_{3}f_1(x)g_3(x))^{2}= 0$. Then there exists $s_4\in S$ such that $s_1s_{2}s_3s_4(f_1(x)g_3(x))=0$, and hence $s_1s_2s_3s_4(f_2(x)g_0(x))=0$. Now obtained equations are given as
		\begin{center}
			\noindent
			$s_{1}(f_0(x)g_1(x))= 0, ~s_{1}(f_1(x)g_0(x))= 0$~~~~~~~~~~~~~~~$(A')$\\
			$s_{2}(f_0(x)g_3(x))=0, ~s_{2}(f_3(x)g_0(x))=0$~~~~~~~~~~~~~~~$(B')$\\
			$\hspace{-1cm}s_1s_2s_{3}s_4(f_2(x)g_0(x)))= 0, s_2s_{3}(f_0(x)g_2(x))= 0$ ~~~~~~~~~~~~$(C')$\\
			$\hspace{3.2cm}s_1s_2s_{3}s_4(f_1(x)g_3(x))=0$.~~~~~~~~~~$(D'')$
		\end{center}
		Now let
		$$f(x)=\sum\limits_{i=0}^{n} \begin{bmatrix}
		a_{i} & b_{i} & c_{i}  \\
		0 & a_{i} & d_{i} \\
		0 & 0 & a_{i} 
		\end{bmatrix}x^{i}, ~and ~g(x)=\sum\limits_{j=0}^{m} \begin{bmatrix}
		a'_{j} & b'_{j} & c'_{j}  \\
		0 & a'_{j} & d'_{j} \\
		0 & 0 & a'_{j} 
		\end{bmatrix}x^{j},$$ where $f_{0}=\sum\limits_{i=0}^{n}a_{i}x^{i}$, $f_{1}=\sum\limits_{i=0}^{n}b_{i}x^{i}$, $f_{2}=\sum\limits_{i=0}^{n}c_{i}x^{i}$, $f_{3}=\sum\limits_{i=0}^{n}d_{i}x^{i}$ $g_{0}=\sum\limits_{j=0}^{m}a'_{j}x^{j}$, $g_{1}=\sum\limits_{j=0}^{m}b'_{j}x^{j}$, $g_{2}=\sum\limits_{j=0}^{m}c'_{j}x^{j}$, and  $g_{2}=\sum\limits_{j=0}^{m}d'_{j}x^{j}$. Take $s=s_1s_2s_3s_4\in S$. Using equations $(A)$, $(A')$, $(B')$, $(C')$ ~and~ $(D'')$, we obtain $sa_{i}a'_{j}=0$, $sa_{i}b'_{j}=0$, $sb_{i}a'_{j}=0$, $sa_{i}c'_{j}=0$, $sb_{i}d'_{j}=0$, $sc_{i}a'_{j}=0$, $sa_{i}d'_{j}=0$, ~and~ $sd_{i}a'_{j}=0$ for all $i,~j$. This implies that $$ \begin{bmatrix}
		s& s& s \\
		0 & s & s\\
		0 & 0 & s
		\end{bmatrix}\begin{bmatrix}
		a_{i} & b_{i} & c_{i}  \\
		0 & a_{i} & d_{i} \\
		0 & 0 & a_{i} 
		\end{bmatrix} \begin{bmatrix}
		a'_{j} & b'_{j} & c'_{j}  \\
		0 & a'_{j} & d'_{j} \\
		0 & 0 & a'_{j} 
		\end{bmatrix}=0$$ for all $i,~j$ and therefore $E$ is an $u$-$S'$-Armendariz rings.
	\end{proof}
\begin{deff}\label{p1}
A ring $R$ is said to be \textit{$u$-$S$-reduced} if for any $a\in R$ such that $a^n = 0$ for some $n \in \mathbb{N}$, there exists fixed element $s\in S$ such that $sa = 0$.
\end{deff}
	In Definition~\ref{p1}, the element $s$ does not depend on the choice of $a \in R$; it depends only on the ring $R$ itself.

	\begin{eg}\label{p2}
	Let $R=\mathbb{Z}_{24}$, and $S=\{\bar{1}, \bar{2}, \bar{4}, \bar{8}, \overline{16}\}$. Evidently, $N^*(R) = \{\bar{6}, \overline{12}, \overline{18} \}$ is the set of all non-zero nilpotent elements of $R$. Let $x\in N(R)$. Then we have $x^{3}=\bar{0}$. Take $s=\bar{4}\in S$. Then $sx=\bar{0}$, and hence $R$ is an $u$-$S$-reduced ring with respect to $\bar{4}\in S$.
\end{eg}
\begin{theorem}\label{arm}
	Let $R$ be an $u$-$S$-reduced with respect to $s\in S$. Then $R$ is an $u$-$S$-Armendariz.
\end{theorem}
\begin{proof}
	 Let $f(x)=a_0+a_1x+\dots+a_nx^n$, ~$g(x)=b_0+b_1x+\dots+b_mx^m\in R[x]$ be such that $f(x)\cdot g(x)=0$. We can assume that $n=m$. Then we have equations: $a_0 b_0 = 0$;  $a_0 b_1 + a_1 b_0 = 0$; $a_0 b_2 + a_1 b_1+ a_2 b_0 = 0$; $a_0 b_3 + a_1 b_2+ a_2 b_1+a_3 b_0 = 0$; $\cdots; a_0 b_n + a_1b_{n-1}+\cdots + a_n b_{0}= 0$. Hence multiplying the second of these equations by $b_0$ yields $b_{0}a_0b_1+b_{0}a_1b_0=0$ implies $(a_1b_0)^{2}=0$, then there exists $s\in S$ such that $sa_1b_0=0$. Now, multiplying third of these equations by $sb_0$ yields $sb_{0}a_0b_2+sb_{0}a_1b_1+sa_2b^{2}_0=0$ implies $(sa_2b_0)^{2}=0$, then there exists $s\in S$ such that $s^{2}a_2b_0=0$. Similarly, multiplying the $i^{th}$ of these equations by $s^{2}b_0$; we get $s^{2}a_{i}b_{0}=0$, for all $4\leq i\leq n$. Thus $s^{2}a_{i}b_{0}=0$, for all $0\leq i\leq n$. Next, multiplying the $(i+1)^{th}$ of the original equations by $s^{2}a_0$; we get $s^{2}a_{0}b_{i}=0$, for all $1\leq i\leq n$. Now, multiplying by $s^{2}$ in the the original equations, reduces to; $s^{2}a_{1}b_{1}=0$; $s^{2}(a_1b_2+a_2b_1)=0$; $\cdots$; $s^{2}(a_1b_{n-1}+a_2b_{n-2}+\cdots+ a_{n-1}b_1)=0$. Then multiplying $(i)^{th}$ reduced equations by $b_1$; we gets; $s^{2}a_{i}b_{1}=0$, for all $2\leq i\leq n$, and if multiplying $(i)^{th}$ of reduced equations by $a_1$; we gets; $s^{2}a_{1}b_{i}=0$, for all $2\leq i\leq n$. Similar repetition yields; $s^{2}a_{i}b_{j}=0$ for all $1\leq i\leq n$, $1\leq j\leq n$, and therefore $R$ is $u$-$S$-Armendariz ring. 
\end{proof}

	\begin{cor}
		Let R be a $u$-$S$-reduced ring. Then $R(+)R$ is a $u$-$S$-Armendariz.
	\end{cor}
	\begin{proof}
		Let $f(x) = (f_0(x), f_1(x))$ and $g(x) = (g_0(x), g_1(x))$ be elements of $(R(+)R)[x]$ satisfying; $f(x) g(x) = 0$. Write;
		$f(x)=\sum_{i=0}^{m}(a_i, u_i)x^i$ and $g(x) = \sum_{j=0}^{n}(b_j, v_j)x^j$ with corresponding representations for $f_{k}(x)$ and $g_{k}(x)$ for $k\in\{0, 1\}$. From the product definition, we have;
		\begin{align*}
		\text{(A)} & \quad f_0(x)g_0(x) = 0. \\
		\text{(B)} & \quad f_0(x)g_1(x) + f_1(x)g_0(x) = 0.
		\end{align*}
		Multiplying equation (B) by $g_0(x)$ and using (A), we get;
		\[
		g_0(x)f_1(x)g_0(x) = 0.
		\]
		This implies;
		\[
		(f_1(x)g_0(x))^2 = 0
		\]
		and since $R[x]$ is an $u$-$S$-reduced, there exists $s\in S$ such that;
		\[
		\text{(C)}\quad s(f_1(x)g_0(x)) = 0.
		\]
		From (B), this implies;
		\[
		\text{(D)} \quad s(f_0(x)g_1(x)) = 0.
		\]
		Combining (A), (C), and (D), and noting that $R$ is $u$-$S$-Armendariz, by Theorem \ref{arm}, we conclude;
		\[
		sa_ib_j = 0, \quad s^{2}u_ib_j = 0, \quad s^{2}a_iv_j = 0 \quad \text{for all } i \text{ and } j.
		\]
		Thus, it follows that; $(s^{2}, s)(a_i, u_i)(b_j, v_j)= (s^{2}, s )(a_ib_j, a_iv_j + u_ib_j) = 0$ for all  $i$ and  $j$, where $(s^{2}, s)\in S(+)R$ is a multiplicative closed set of $R(+)R$.
	\end{proof}
	Following from \cite{ab22}, a ring $R$ is said to be strongly Hopfian if the chain of annihilators $\text{ann}(a) \subseteq \text{ann}(a^2)\cdots \subseteq\text{ann}(a^{n}) \subseteq\text{ann}(a^{n+1}) \subseteq \cdots$ stabilizes for each $a \in R$. Later, in 2021, Ahmed et al. \cite{ha21} generalized this concept to \( S \)-strongly Hopfian rings as follows:
	
	\begin{deff}$\cite{ha21}$\label{strongly}
		A ring $R$ is called $S$-\textit{strongly Hopfian ring} if the chain of annihilators $\text{ann}(a) \subseteq \text{ann}(a^2) \subseteq \cdots\subseteq\text{ann}(a^n)\subseteq\text{ann}(a^{n+1}) \subseteq \cdots$ is $S$-stationary for each $a \in R$.	
	\end{deff}
	\begin{eg}
		$S$-integral domains and $S$-Noetherian rings are $S$-strongly Hopfian.
	\end{eg}

	\begin{prop}\label{str}
		Each $S$-reduced ring is $S$-strongly Hopfian.
	\end{prop}
	\begin{proof}
		Let $R$ be an $S$-reduced ring and $a\in R$. Now consider the increasing chain of annihilators; $$\text{ann}(a) \subseteq \text{ann}(a^2) \subseteq \cdots\subseteq\text{ann}(a^n)\subseteq\text{ann}(a^{n+1}) \subseteq \cdots,$$ where $n\geq 1$. Let $y\in \text{ann}(a^{n+1})$. Then $ya^{n+1}=0$, and $(ya)^{n+1}=0$. Consequently, there exists $s\in S$ such that $sya=0$, and $sy\in ann(a)$. Thus $s(ann(a^{n+1}))\subseteq ann(a^{n})$ for all $n\geq 1$.
	\end{proof}
	\begin{cor}
		For each ring $R$, the ring $R/Nil_{S}(R)$ is $S$-strongly Hopfian.
	\end{cor}
	\begin{deff}$\cite{ba24}$\label{pure}
		An ideal $I$ of $R$ is called $S$-\textit{pure} if for all $a\in I$, there
		exist $b\in I$ and $s\in S$ such that $sa = ab$.
	\end{deff}
	Notice that every pure ideal is an $S$-pure, but the converse need not be true (see \cite[Example 3.3]{ba24}).
	
%
	\begin{deff}
		A ring \( R \) is called an \( S \)-\textit{PF ring} if, for every element \( a \in R \), the annihilator \( (0:a)_R \) is an \( S \)-pure ideal. This means that for each \( x \in (0:a) \), there exist an element \( y \in (0:a) \) and \( s \in S \) such that \( xy = sx \).
	\end{deff}
	\begin{eg}\label{pf}
		Every $S$-PF-ring is an $S$-reduced ring.
	\end{eg}
	\begin{proof}
		Let $R$ be an $S$-PF-ring and $a\in R$ with $a^{n}=0$ for some $n\in\mathbb{N}$. Then $a\in(0:a^{n-1})$ and $(0:a^{n-1})$ is $S$-pure, there exist $b\in (0:a^{n-1})$ and $s_{1}\in S$ such that $ab=s_{1}a$. As $a^{n-1}b=0$, $aba^{n-2}=s_{1}a^{n-1}=0$. Consequently, $a\in(0:s_{1}a^{n-2})$. Since $(0:s_{1}a^{n-2})$ is $S$-pure, there exist $b_{1}\in (0:s_{1}a^{n-2})$ and $s_{2}\in S$ such that $ab_{1}=s_{2}a$. As $s_{1}a^{n-2}b_{1}=0$, $s_{1}s_{2}a^{n-2}=0$ implies that $a\in(0:s_{1}s_{2}a^{n-3})$. Using similar process, we get $a\in(0:s_{1}s_{2}\ldots s_{n-3}a^{(n-\{n-2\})}))$ for all $n\geq 3$, $s_{1}s_{2}\ldots s_{n-3}a=0$. Put $t=\prod\limits_{k=1}^{n-3}s_{k}$. Thus $ta=0$, as desired.
	\end{proof}
	Recall from \cite[Definition 1.6.10 ]{fw16}, an $R$-module $M$ is called an $S$-\textit{torsion module} if for any $m\in M$, there exists $s\in S$ such that $sm=0$. Following from \cite{ki24}, an $S$-prime ideal $P$ of $R$ is said to be an $S$-\textit{minimal} $S$-\textit{prime} if whenever $Q\subseteq P$ for some $S$-prime ideal $Q$ of $R$, then $sP\subseteq Q$ for some $s\in S$. By \cite[Proposition 2.1]{ki24}, $Spec_{S}R$ contains a prime ideal which is disjoint from $S$, and every prime ideal contains a minimal prime ideal. Also every minimal prime ideal of $R$ disjoint from $S$ is $S$-minimal $S$-prime, by \cite[Example 2.6]{ki24}. Then the intersection of all $S$-prime ideals of $R$ is equal to the intersection of all $S$-minimal $S$-prime ideals of $R$. If $R$ is an $S$-reduced, then the intersection of all $S$-minimal $S$-prime ideals of $R$ is $S$-zero.

	\begin{rem}$\cite{wq23}$
		Let $M$ and $N$ be $R$-modules. Then we have the following:
		\begin{enumerate}
			\item An $R$-homomorphism $f:M\longrightarrow N$ is an $S$-injective if and only if $Ker(f)$ is an $S$-torsion.
			\item An $R$-homomorphism $f:M\longrightarrow N$ is an $S$-surjective if and only if $Coker(f)$ is an $S$-torsion.
		\end{enumerate}
	\end{rem}
	\begin{deff}
		Let $\wedge$ be an indexed set, and $R_{i\in\wedge}$ be rings. We define the natural projection for each $j \in \wedge$ as $\pi_j: \prod_{i\in\wedge} R_i \longrightarrow R_j$. Then we say that $R$ is $S_{i}$-subdirect product of $R_i$, where $S_{i}$ be multiplicative closed subset of $R_{i}$ for all $i\in\wedge$, if the following statements are satisfied:
		\begin{enumerate}
			\item There exists an $S$-injective ring homomorphism $f: R \longrightarrow \prod_{i\in\wedge} R_i$,
			\item  For every $j \in\wedge$ the map $\pi_jof: R \longrightarrow R_j$ is $S$-surjective.
		\end{enumerate}
	\end{deff}
	
	\begin{theorem}$($\textbf{Structure theorem For $S$-reduced rings}$)$\label{struct} Let $\{P_{i}\}_{i\in\wedge}$ be $S$-minimal $S$-prime ideals of $R$. Then R is an $S$-reduced if and only if $R$ is an $S$-subdirect product of   $\overline{S}_{i}$-integral domains, where $\overline{S}_{i}=\{s+P_{i}\mid s\in S\}$ is a multiplicative closed subset of $R/P_{i}$.
	\end{theorem}
	\begin{proof}
		Consider a ring homomorphism $f: R\longrightarrow \prod_{i\in\wedge} R/P_i$ defined as $f(r)=(r+P_1, ~r+P_2, \ldots)$, where $P_{i}\in Spec_{S}R$ for all ${i\in\wedge}$. Then $Kerf=\{r\in R\mid f(r)=P_{i}\}=\{r\in R\mid r\in P_{i}  ~for ~all ~i\in\wedge\}=\{r\in R\mid r\in\bigcap_{i\in\wedge} P_{i}=\bigcap_{P\in Spec_{S}(R)} P\}$. Since $R$ is an $S$-reduced, by Theorem \ref{szero}, $Kerf=\{r\in R\mid sr=0 ~for ~some ~s\in S\}$. Thus $kerf$ is $S$-torsion, and hence $f$ is $S$-injective. For every $j\in\wedge$, we define a ring homomorphism $\pi_{j}o f:R\longrightarrow R/P_{j}$, where $\pi_{j}:\prod_{i\in\wedge} R/P_i\longrightarrow R/P_{j}$ is a natural projection. Evidently, for all $j\in\wedge$,  $Im(\pi_{j}o f)=R/P_{j}$, for this, let $r\in R$, $(\pi_{j}o f)r=\pi_{j}(r+P_{1}, r+P_{1}, \ldots)=r+P_j$ for all $r\in R$. Thus $Coker(\pi_{j}o f)$ is an $S$-torsion. Therefore $R$ is $S$-subdirect product of $R/P_{i}$. Since $R/P_{i}$ is $\overline{S}_{i}$-integral domain for all $i\in\wedge$, and hence $R$ is $S$-subdirect product $\overline{S}_{i}$-integral domain.
		
		Conversely suppose $R$ is $S$-subdirect product of $R/P_{i}$, i.e,  $\overline{S}_{i}$-integral domain for all $i\in\wedge$, and $f: R \longrightarrow \prod_{i\in\wedge} R/P_{i}$ defined as $f(r)=(r+P_1, ~r+P_2, \ldots)$ is an $S$-injective ring homomorphism. Let $r\in R$ and $r^n=0$ for some $n\in\mathbb{N}$. Then we have $f(r^{n})=f(0)$, $(f(r))^{n}=(P_{1}, P_{2},\ldots)$. Consequently, $r^n\in P_{i}$ for all $i\in\wedge$, and so $(s+P_{i})(r+P_{i}) = P_{i}$ for some $s\in S$ because every $R/P_{i}$ is an $\overline{S}_{i}$-integral domain. Hence $sr\in P_{i}$ for all $i\in\wedge$, and so $sr\in\bigcap_{i\in\wedge} P_{i}$ implies $sr\in Kerf$. Then there exists $s'\in S$ such that $ss'r=0$. Thus $R$ is an $S$-reduced.
	\end{proof}

	\noindent
	\textbf{Conflict of Interest Statement:} The authors declare that they have no conflict of interest.\\

\end{document}